\documentclass{article}
\usepackage{arxiv}
\usepackage{amsthm} 
\usepackage{amssymb}
\usepackage{amsmath}
\usepackage{empheq}
\usepackage{amscd}
\usepackage{graphicx}
\usepackage{graphics}
\usepackage[noadjust]{cite}
\usepackage{amsthm}
\usepackage{caption}
\usepackage{multirow}
\usepackage{array}
\usepackage{rotating}
\usepackage{hyperref}

\usepackage{indentfirst}
\usepackage{tikz}
\usepackage{calc}
\usepackage{epsfig}
\usepackage[numbers]{natbib}
\usepackage[makeroom]{cancel}
\usepackage{multicol}
\usepackage{csquotes}
\usepackage{enumitem}

\usepackage{hyperref}       
\hypersetup{colorlinks,linkcolor={green},citecolor={green},urlcolor={black}}
\usepackage{optidef}

\renewcommand{\bold}[1]{#1}
\renewcommand{\boldsymbol}[1]{#1}

\newcommand{\ur}{\bold{u}}
\newcommand{\up}{\bold{u}}
\newcommand{\um}{\bold{u}}
\renewcommand{\u}{\bold{u}}

\newcommand{\LamAL}{{\bold{w}}}
\newcommand{\Lamm}{\bold{v}}
\newcommand{\LamALn}{{\boldsymbol{\varepsilon}}}
\newcommand{\Lam}{{\boldsymbol{\lambda}}}

\newcommand{\scs}{\boldsymbol{\varphi}}

\newcommand{\regp}{\mu}

\newcommand{\eref}[1]{eq. \eqref{#1}}                  
\DeclareMathOperator*{\minimize}{minimize}

\DeclareMathOperator*{\minimax}{minimax}

\begin{document}
\title{Full Waveform Inversion and Lagrange Multipliers}

\author{\href{https://orcid.org/0000-0002-9879-2944}{\hspace{1mm}Ali Gholami} \\
  Institute of Geophysics, Polish Academy of Sciences, Warsaw, Poland\\
  \texttt{agholami@igf.edu.pl} \\ 
  }

\graphicspath{{./figs/}}

\renewcommand{\shorttitle}{FWI and Lagrange Multipliers ~~~~~~~~~~~~~~~~~~~~~~~~~ Ali Gholami}

\maketitle
\begin{abstract}
Full-waveform inversion (FWI) is an effective method for imaging subsurface properties using sparsely recorded data. It involves solving a wave propagation problem to estimate model parameters that accurately reproduce the data. 
Recent trends in FWI have led to the development of extended methodologies, among which source extension methods leveraging reconstructed wavefields to solve penalty or augmented Lagrangian (AL) formulations have emerged as robust algorithms, even for inaccurate initial models. Despite their demonstrated robustness, challenges remain, such as the lack of a clear physical interpretation, difficulty in comparison, and reliance on difficult-to-compute least squares (LS) wavefields.
This paper is divided into two critical parts. In the first, a novel formulation of these methods is explored within a unified Lagrangian framework.  This novel perspective permits the introduction of alternative algorithms that employ LS multipliers instead of wavefields. These multiplier-oriented variants appear as regularizations of the standard FWI, are adaptable to the time domain, offer tangible physical interpretations, and foster enhanced convergence efficiency. The second part of the paper delves into understanding the underlying mechanisms of these techniques. This is achieved by solving the FWI equations using iterative linearization and inverse scattering methods. The paper provides insight into the role and significance of Lagrange multipliers in enhancing the linearization of FWI equations. It explains how different methods estimate multipliers or make approximations to increase computing efficiency. Additionally, it presents a new physical understanding of the Lagrange multiplier used in the AL method, highlighting how important it is for improving algorithm performance when compared to penalty methods.
\end{abstract}


%

\section{Introduction}
Full waveform inversion (FWI) is a widely employed and effective technique for imaging subsurface properties in geophysical studies \citep{Tarantola_1984_ISR,Pratt_1998_GNF,Virieux_2009_OFW,Operto_2023_FWI}. This involves solving a wave propagation problem to determine the model parameters that accurately reproduce the recorded/observed data. The governing law is defined by a set of nonlinear equations that capture the wavefields and establish their connection to the model parameters, thereby enabling the calculation of synthetic data. The primary objective of FWI is to solve this set of nonlinear equations to estimate the model parameters based on observed data.

Unlike linear systems, such nonlinear systems of equations cannot be solved by direct methods due to their complexity. Therefore, they must be solved using iterative methods. These iterative methods can be generated based on various principles, such as iterative linearization approaches or inverse problem theory. In the iterative linearization approach, the nonlinear system is linearized at each iteration, and the linearized equations are then solved to obtain an updated solution. This process is repeated until convergence is achieved. On the other hand, the inverse problem theory approach replaces the problem of solving the system of nonlinear equations with the optimization of some nonlinear functional. The desired solution is then obtained by finding the best solution in a precise sense, usually by minimizing the difference between the observed data and the predicted data \citep{Tarantola_2005_IPT}. 
This paper is divided into two independent sections, each of which focuses on one of these frameworks of solving the FWI equations.

In a first section, we consider the traditional form of FWI in the framework of inverse problem theory. It addresses the estimation of subsurface model parameters by minimizing the difference between observed and synthetic data through a nonlinearly constrained least-squares (LS) optimization problem. To incorporate physical constraints and prior information for regularization, the widely adopted Lagrange multiplier method is employed \citep{Haber_2000_OTS}. This approach introduces Lagrange multipliers to enforce constraints such as the wave equation governing the wave propagation problem. The desired solution, which encompasses the model parameters, wavefields, and multipliers, is obtained by solving three sets of nonlinear equations that satisfy the Karush-Kuhn-Tucker (KKT) optimality conditions \citep{Nocedal_2006_NO}. These conditions are satisfied by setting the partial derivatives of the Lagrangian with respect to the multipliers, wavefields, and model parameters equal to zero.

The solution of the KKT equations can be achieved through Newton-type methods. These methods iteratively linearize the system around the current solution estimate and solve the resulting linear system using preconditioned Krylov iterative methods \citep{Akcelik_2002_MNK, Hoffmann_2021_PFA}. However, efficiently solving large KKT systems poses significant challenges from both computational and memory perspectives. The large size of the systems can lead to high computational costs and memory requirements, especially for large-scale problems encountered in practical applications. Furthermore, the Hessian matrix associated with the KKT system may be ill-conditioned and indefinite, leading to numerical instabilities and convergence issues. To address the ill-conditioning, appropriate preconditioning or regularization procedures are needed, as well as effective ways for memory issues.

The memory issue can be mitigated through variable projection \citep{Golub_2013_VPM}. Traditionally, this approach involves eliminating the wavefields and multipliers from the optimization variables by exactly satisfying two sets of the KKT equations, the forward and adjoint wave equations. The gradient of the reduced objective function, necessary for updating the model parameters, is then computed by correlating the wavefields and multipliers, which propagate forward and backward in time, respectively \citep{Plessix_2006_RAS}. While this reduced approach offers significant computational and memory savings, it still faces challenges related to the ill-conditioning of the inversion process, leading to numerical instabilities and slow convergence \citep{Metivier_2017_TRU}.
Recent advancements in computational resources have enabled solving large and complex FWI problems using advanced numerical optimization methods. As a result, addressing the ill-conditioning issue and improving the convergence behavior of FWI methods have become an active areas of research \citep[e.g.,][]{VanLeeuwen_2013_MLM, Yang_2017_AOT, Warner_2016_AWI, Li_2016_FWI,Huang_2018_VSE,Aghamiry_2019_IWR, Alkhalifah_2019_AEW,Rizzuti_2021_DFW, Gholami_2022_EFW,Operto_2023_FWI,Barnier_2023_FWIt,Metivier_2022_ARO}.

This study focuses on investigating FWI methods within a Lagrangian framework, including the penalty and augmented Lagrangian (AL) methods. In its standard form, the penalty method enforces constraints by adding a quadratic term to the objective function that penalizes the sum of the squares of the wave-equation violations. The amount of penalization is determined by the penalty parameter, which must be increased with iteration to satisfy the constraints. To solve the penalty objective function, \citet{VanLeeuwen_2013_MLM} employed the Wavefield Reconstruction Inversion (WRI) algorithm, which follows a two-step alternating iteration for updating the wavefields and model parameters.
The main issue with WRI is that it requires high values of the penalty parameter to satisfy the wave equation, yet high values slow down the algorithm's convergence because this method essentially employs a step length that varies with the inverse of the penalty parameter.
To solve this issue, the AL method combines both the Lagrangian and penalty methods. This results in the iteratively refined WRI (IR-WRI) algorithm, similar to WRI but with a modified volume source term \citep{Aghamiry_2019_IWR}.
The AL method satisfies the wave equation even for moderate values of the penalty parameter without the need for excessively large values. This preserves the well-posedness of the algorithm through iterations which leads to accurate estimations of the multiplier and to an improved convergence behavior.

In recent years, numerous papers have demonstrated the superior performance of Lagrangian-based methods, such as WRI and IR-WRI, compared to standard FWI, using various numerical examples \citep[e.g.][]{VanLeeuwen_2013_MLM,Huang_2018_VSE, Aghamiry_2019_IWR, Alkhalifah_2019_AEW, Rizzuti_2021_DFW, Gholami_2022_EFW, Operto_2023_FWI, Lin_2023_FWR}. However, despite their effectiveness, these methods face several challenges:
(1) They are based on least squares (LS) wavefields, also known as extended wavefields or data-assimilated wavefields. These wavefields are obtained by solving an augmented wave equation that simultaneously incorporates the wave equation and the data equation in the LS sense. Efficient computation of these LS wavefields in time-domain applications is very challenging because explicit time-stepping methods for solving the augmented wave equation are lacking \citep[see][and references therein]{Wang_2016_FIR,Rizzuti_2021_DFW,Gholami_2022_EFW,Operto_2023_FWI,Lin_2023_FWR}
(2) Unlike traditional FWI methods, these approaches do not involve a backward wavefield. However, the model update still relies on correlating the LS wavefields and a quantity that is a function of these wavefields whose physical meaning is not immediately clear.
(3) Interpreting the resulting algorithms physically can be challenging, making the mechanism of each method unclear or sometimes misunderstood. This lack of clarity makes it difficult to compare WRI and IR-WRI with each other and with standard FWI. This issue is mainly because they have not been considered under a unified formulation, and the physical significance of Lagrange multipliers and their function in these methods are not well understood.

To address these challenges, in this study, by investigating these methods within a unified framework, we aim to shed light on the underlying mechanisms and explore their potential for more practical applications. Understanding the physical meaning of Lagrange multipliers and their role in the inversion process not only enhances our theoretical understanding of the methods but also opens up new avenues for improving their performance and applicability in practice.
We demonstrate that there exist alternative formulation for WRI and IR-WRI which gives the same solution to the inverse problem, but they are based on the computation of LS multipliers rather than LS wavefields. 
They lead to inversion algorithms that are extremely similar to the standard FWI, can be performed in the time domain using explicit time-stepping methods, and are physically interpretable.
To achieve this, we formulate the FWI problem as a two-step iterative process. The first step involves solving a convex-concave min-max optimization problem to compute the wavefields and multipliers. Subsequently, the model parameters are updated in the second step using a correlation of the computed wavefields and multipliers, akin to the standard FWI. Solving the first step requires handling a two-by-two block linear system of equations, also known as a saddle point system. The penalty and AL methods are two different approaches for regularization of the ordinary Lagrangian hence enable the stable solution of this saddle point system. For solving this system, two distinct approaches can be employed: the wavefield-oriented approach and the multiplier-oriented approach. These approaches emphasize either wavefield reconstruction or multiplier computation as the primary goal, providing flexibility in implementing the algorithms.

The wavefield-oriented approach, which reduces to WRI in the penalty method \citep{VanLeeuwen_2013_MLM} and to IR-WRI in the AL method \citep{Aghamiry_2019_IWR}, focuses on reconstructing the LS wavefields as an initial step. This approach aims to obtain accurate wavefield information, which is then used to derive the Lagrange multipliers. By utilizing the LS wavefields, multipliers can be estimated cheaply and used to calculate the model update.
By contrast, the multiplier-oriented approach prioritizes the calculation of LS multipliers as the primary step. It utilizes these multipliers as secondary volume sources in the wave equation to calculate the wavefields. The advantage of this approach is that it allows the use of standard time-stepping methods to compute the wavefields and multipliers.

In the second part of this paper, we examine how to solve the original nonlinear equations that govern FWI without using optimization theory. Instead, we solved the FWI equations using scattering concepts and iterative linearization methods, aiming to find a set of model parameters and wavefields that satisfy the wave equation while remaining consistent with the recorded data.
In this framework, Lagrange multipliers are introduced as auxiliary variables that split the nonlinear equations, hence simplifying the linearization process.
This allows us to shed light on the physical interpretation of Lagrange multipliers, revealing their role as scattering sources or scattered data, which significantly contributes to the linearization of FWI equations and improves the overall iteration efficiency. We demonstrate how various FWI methods strive to estimate accurate Lagrange multipliers or employ approximations to enhance computational and memory efficiency. Importantly, this study offers a pioneering insight into the physical meaning of the Lagrange multiplier used in the AL method, emphasizing its significance in enhancing algorithm performance compared to penalty method. 


\section{NOTATION}
In this paper, we develop iterative methods for solving FWI. We use a streamlined notation to represent the iterative process to maintain a clear and concise presentation throughout the analysis. In particular, we denote the background values of the variables without indices, while the updated variables are denoted by a superscript ``$+$". In some circumstances, it may be necessary to construct updated variables using the same variables used to create background variables. The superscript ``$-$" is used to identify these variables. For instance, we input the variables' background values ($\bold{m}$ and $\bold{u}$), as well as their previous values ($\bold{m}^-$ and $\bold{u}^-$) for each iteration. The updated values of the variables, denoted as $\bold{m}^+$ and $\bold{u}^+$, are obtained after performing the necessary computations and updates. The updated variables from the previous iteration serve as the background variables for the subsequent iteration, which adheres to the same convention. 
We hope to simplify the presentation and make the iterative process easier to understand by using this notation. It keeps the text and formulas simple, and allows us to concentrate on the key components of each iteration. 

\section{Lagrangian based full waveform inversion}
We start with the assumption that the seismic wavefield, denoted by $u(t,\bold{x};\bold{x}_s)$, is generated by a source term $b(t,\bold{x};\bold{x}_s)$ in a medium characterized by variable velocity and constant density. This wavefield satisfies the wave equation given by
\begin{equation} \label{waveEqc}
m(\bold{x}) \frac{\partial^2}{\partial t^2}u(t,\bold{x};\bold{x}_s) - \nabla^2 u(t,\bold{x};\bold{x}_s)=b(t,\bold{x};\bold{x}_s),
\end{equation}
where $m(\bold{x})$ represents the model parameters, specifically, the squared slowness, which varies with the spatial coordinate $\bold{x}$ in two or three dimensions. Source locations are denoted by $\bold{x}_s$ for $s=1,2,...,N_s$, where $N_s$ is the total number of sources. Operator $\nabla^2$ represents the Laplacian.

In practical applications, the wavefield is typically sampled at the receiver locations $\bold{x}_r$ for $r=1,2,...,N_r$ (where $N_r$ is the number of receivers). These receivers can be positioned in various ways, depending on the imaging objectives and characteristics of the subsurface under investigation. Examples include placing receivers in boreholes, deploying them on the seafloor, and distributing them throughout the survey area to capture the seismic wavefield as it propagates through the subsurface. The recorded samples at these receiver locations, denoted by $d(t,\bold{x}_r;\bold{x}_s)$, are equivalent to $u(t,\bold{x}_r;\bold{x}_s)$ and serve as the input data for the FWI algorithms to estimate the model parameters.

For simplicity, our analysis is based on single-source problems. However, this does not impose any limitations on our formulation because the extension to multiple sources can be achieved by summing the final model update formula over the sources. 

Assuming that a known source $\bold{b}^*$ is used to generate the true wavefield $\bold{u}^*$ in a medium characterized by model parameters $\bold{m}^*$, we can record the corresponding data $\bold{d}^*$. These quantities are connected by a pair of coupled equations.
\begin{subequations} \label{Pair1}
\begin{align}
&\bold{A}(\bold{m}^*)\bold{u}^* =\bold{b}^*, \label{waveEq}\\
&\bold{P}\bold{u}^*=\bold{d}^*. \label{dataEq}
\end{align}
\end{subequations}
Here, $\bold{A}(\bold{m}) \in \mathbb{R}^{N \times N}$ represents the discretization of wave equation \eqref{waveEqc} with appropriate boundary conditions. The matrix size $N$ is determined by $N_m \times N_t$, where $N_m$ is the number of grid points used to sample the subsurface model and $N_t$ is the number of time samples. Meanwhile, $\bold{P} \in \mathbb{R}^{M \times N}, M=N_r\times N_t,$ represents the sampling operator, which samples the wavefield at the receiver locations $\bold{x}_r$.

To solve nonlinear equations \eqref{waveEq} and \eqref{dataEq}, they are often transformed into a nonlinear LS problem using inverse problem theory. Two popular approaches were considered.
The first approach considers the model parameters $\bold{m}^*$ as the only independent variable and considers the wavefield $\bold{u}^*$ as a function of $\bold{m}^*$ by exactly satisfying the wave equation, that is, $\bold{u}^* = \bold{A}(\bold{m}^*)^{-1}\bold{b}^*$ \citep{Tarantola_1984_ISR,Gauthier_1986_TDN,Pratt_1998_GNF}. This reduced-space approach has the advantage of a relatively small problem size because the model parameters are independent of source and time. However, the associated formulation can be ill-conditioned.
The interested reader is referred to \citep{Gholami_2023_MWI} for AL-based formulation of the approach.
The second approach considers both the model parameters and the wavefield as independent optimization variables. This formulation does not require precise satisfaction of the wave equation at each iteration of the optimization process. Instead, it only needs to be satisfied in the final solution \citep{VanLeeuwen_2013_MLM,Aghamiry_2019_IWR,Rizzuti_2021_DFW,Gholami_2022_EFW}. This formulation can act as a form of regularization, reducing the ill-conditioning of the optimization problem. However, because the wave equation is only approximately solved in each iteration, it imposes additional computational and memory burdens, which are significant disadvantages.


Following the second approach, we can estimate the model and wavefield variables simultaneously from the data by solving
\begin{equation} \label{EFWI_obj}
\minimize_{\bold{m},\bold{u}}~\|\bold{Pu}-\bold{d}^*\|_2^2\quad 
\text{subject to}\quad \bold{A(m)u}=\bold{b}^*.
\end{equation}
Indeed, in practical applications of FWI, regularization terms are often included in the objective function to ensure stable model updates and to impose constraints on the updated model parameters. These regularization terms help filter out undesired features and restrict the model to desired spaces \citep{Gholami_2023_MWI}.
However, for the purpose of interpretation and to maintain the focus on the role of the multipliers, we have chosen not to include such regularization terms in our formulas. This allows us to develop and present iterative methods in their most basic form, making it easier to analyze the contributions of the various components involved in iterations.
However, it is important to note that the iterative methods presented here can easily be extended to include regularization terms. Regularization terms would introduce additional terms into the model update equations. Despite the addition of regularization, the core structure and principles of iterative methods remain unchanged, and the physical meaning of the multipliers remains important for understanding the mechanisms.

The Lagrange multiplier method formulates \eref{EFWI_obj} as
\begin{equation}  \label{L}
\min_{\bold{m},\u}\max_{\Lamm} \mathcal{L}(\bold{m},\u,\Lamm)=\frac12\|\bold{P}\u-\bold{d}^*\|_2^2+\langle\Lamm,\!\bold{A(m)}\u\!-\bold{b}^*\rangle_{t,x}
\end{equation}
where $\Lamm\in \mathbb{R}^{N}$ denotes the Lagrange multiplier vector and $\langle\cdot,\cdot\rangle_{t,x}$ is the canonical inner product with the integration over space and time.  The KKT conditions \citep{Nocedal_2006_NO}, which are first-order necessary conditions  for optimality of the solution, require that the partial derivatives of the Lagrangian with respect to $\Lamm$, $\u$, and $\bold{m}$ be zero at the optimum point, which is a saddle point. 
Three sets of nonlinear equations are produced and must be solved.
Newton-type methods are commonly employed to find the optimum solution, where the associated KKT system is solved simultaneously for all variables \citep{Haber_2000_OTS,Akcelik_2002_MNK,Hoffmann_2021_PFA}. In other words, a Krylov iterative method is used to solve the whole KKT system during each Newton iteration. The challenge with this method is that the Hessian matrix is very large and might be ill-conditioned and indefinite, which in turn requires appropriate preconditioning. Alternatively, this paper introduces an iterative process that breaks down the optimization problem into a saddle point (min-max) subproblem involving $\Lamm$ and $\u$. Subsequently, a simple preconditioned gradient descent step is performed to update the model parameters.
The resulting iteration is given by:
\begin{subequations}
\begin{align} 
&(\u^{+},\Lamm^{+})=\arg \minimax_{\u,\Lamm} \mathcal{L}(\bold{m},\u,\Lamm), \label{fwi_min_max}\\
&{m}^+={m}-\alpha \frac{\langle \Lamm^+, \partial_{tt} \u^{+}\rangle_t }{\langle \partial_{tt} \u^{+}, \partial_{tt} \u^{+}\rangle_t} .\label{fwi_dm}
\end{align}
\end{subequations}
Here, $\alpha>0$ represents the step length. The inner product $\langle \cdot, \cdot \rangle_t$ is performed only over time. The expression $\langle \Lamm^+, \partial_{tt} \u^{+}\rangle_t$ denotes multiplication and summation over time of the multiplier and the second time-derivative wavefield, which serves as the gradient of the Lagrangian function with respect to the model parameters. The denominator in the rightmost term in \eref{fwi_dm} acts as a preconditioning factor, and the division is performed element-wise.

Equation \eqref{fwi_dm} indeed offers valuable insights into the structure of the multiplier. To ensure the correct dimensionality for the model update, the multiplier should indeed involve a multiplication of the wavefield and a model perturbation. More specifically, a model update $\delta\bold{m}=\bold{m}^+-\bold{m}$ obtained from the evaluation of the division term in \eref{fwi_dm} shows that the multiplier should be of the form $\Lamm^{+}=-\delta\bold{m}\circ \partial_{tt} \u^{+}$, where $\circ$ denotes the element-wise product in space. This form of the multiplier represents a scattering source or Born secondary source \citep{Tarantola_1988_TBI}. In Section \ref{physical_interp} of this paper, we investigate this physical meaning of the multipliers.

Once the updated wavefield and multiplier are available, calculating the updated model in \eref{fwi_dm} is straightforward. 
Hence, all the efforts are toward developing robust methods for estimating wavefield and multiplier during the iterative process. In the following subsections, we will explore different FWI algorithms for the estimation of $\u^{+}$ and $\Lamm^{+}$.

\subsection{The Standard Lagrangian}
We may use the original Lagrangian to calculate the updated wavefield and multiplier. The necessary conditions for optimality of them are
\begin{equation}
\frac{\partial \mathcal{L}(\bold{m},\u,\Lamm)}{\partial \Lamm} = \bold{0} \quad \text{and} \quad
\frac{\partial \mathcal{L}(\bold{m},\u,\Lamm)}{\partial \u} = \bold{0},
\end{equation}
which result in the following two-by-two block linear system, so called saddle point system,
\begin{equation} \label{redused}
\begin{pmatrix}
\bold{A}(\bold{m}) & 0 \\
\bold{P}^T\bold{P} & \bold{A}(\bold{m})^T
\end{pmatrix}
\begin{pmatrix}
\ur^{+} \\
\Lamm^{+}
\end{pmatrix}=
\begin{pmatrix}
\bold{b}^* \\
\bold{P}^T\bold{d}^*
\end{pmatrix}.
\end{equation}
The first equation is independent of the multiplier, which allows the wavefield to be calculated by solving the wave equation. The multiplier is then obtained by solving the adjoint equation.  This results in the following iteration: 
\begin{subequations} \label{redused2}
\begin{align} 
&\ur^{+}=\bold{A}(\bold{m})^{-1}\bold{b}^*,\label{redused2a}\\
&\Lamm^{+} = \bold{A}(\bold{m})^{-T}\bold{P}^T(\bold{d}^*-\bold{P}\ur^{+}), \label{redused2b}\\
&{m}^+={m}-\alpha \frac{\langle \Lamm^+, \partial_{tt} \u^{+}\rangle_t }{\langle \partial_{tt} \u^{+}, \partial_{tt} \u^{+}\rangle_t} . \label{redused2c}
\end{align}
\end{subequations} 
These equations define a wavefield and a multiplier that can be used as implicit functions of $\bold{m}$ in the Lagrangian function, \eref{L}, reducing it to a function involving only $\bold{m}$ as in \eref{FWI_obj}.
It can be shown that $\langle \Lamm^+, \partial_{tt} \u^{+}\rangle_t$ is indeed the gradient of the reduced function in \eref{FWI_obj} and thus \eref{redused2c} is a preconditioned gradient descent step for its minimization \citep{Tarantola_1984_ISR,Gauthier_1986_TNI,Plessix_2006_RAS,Virieux_2009_OFW}.
This implies that, in the standard FWI, the first two sets of KKT conditions are automatically satisfied. The main objective of the subsequent iterations is to satisfy the final condition, which is highly challenging owing to the poor conditioning of the problem. 

Regarding the estimated multipliers, it is crucial to acknowledge that \eref{redused2b} provides only an approximate multiplier estimate with incorrect amplitudes, which can be out of range unless data residuals are small. This discrepancy arises from the fact that the modeling operator $\bold{P}\bold{A}(\bold{m})^{-1}$ is not unitary, and therefore, the estimated multiplier in \eref{redused2b} may not precisely fit the scattered data. In practice, this issue can impact the convergence and accuracy of the inversion process. However, we can address this challenge by tuning the step length $\alpha$ appropriately or using prior information to enforce the desired bounding of the multipliers within a certain range. Additionally, we can modify the Lagrangian formulation to include damping terms that penalize large Lagrange multiplier values, which helps stabilize the optimization process. 

%
In the following subsections, we introduce the penalty and AL methods as means of regularizing the Lagrangian optimization process. However, our formulation differs slightly from the standard penalty and AL approaches, aiming to offer a simpler and more interpretable formulation while being more general. Notably, we retain the original Lagrange multiplier $\Lamm$ in our formulation, in contrast to the conventional practice of eliminating it to obtain standard forms. It is important to note that our formulation do not introduce any additional computational or memory burdens to the final algorithms. On the contrary, it allows for a more intuitive interpretation and understanding of the optimization process.


\subsection{Regularization by the quadratic penalty method}
The quadratic penalty method can be seen as introducing a damping term in the Lagrangian function to penalize large multiplier values \citep{Gill_2012_PDA,Wright_2022_ODA}. This leads to a modified formulation of the min-max problem in \eref{fwi_min_max}:
\begin{equation}  \label{LP}
(\u^{+},\Lamm^{+})=\arg \minimax_{\u,\Lamm} ~ \mathcal{L}(\bold{m},\u,\Lamm) - \frac{1}{2\regp}\|\Lamm\|_2^2.
\end{equation} 
Here, $\regp>0$ controls the strength of the damping applied to the multiplier.
This is simply zeroth-order Tikhonov regularization of the maximization problem. The minus sign of the regularization term comes from the maximization of the Lagrangian with respect to the multiplier.
As $\regp$ tends towards infinity, the objective function approaches the original Lagrangian function, maintaining its properties. Conversely, as $\regp$ approaches zero, the multipliers dampen towards zero.
By adjusting $\regp$, we can fine-tune the influence of the damping term and enforce the desired bounding of the multipliers within a certain range. 

This approach is closely connected to the standard quadratic penalty method, but it offers a simpler and more interpretable formulation.
Note that in the maximization step, the Lagrange multiplier $\Lamm$ can be explicitly computed as $\Lamm = \regp (\bold{A}(\bold{m})\ur - \bold{b}^*)$. Substituting this expression into the objective function \eqref{LP} eliminates the Lagrange multiplier, simplifying the formulation to a standard quadratic penalty function in \eref{Penalty_obj}.
Also, we can eliminate $\u$ from this penalty function yielding a function solely dependent on $\bold{m}$ \citep{vanLeeuwen_2016_PMP,vanLeeuwen_2019_ANO,Symes_2020_WRI, Rizzuti_2021_DFW, Gholami_2022_EFW}. 

The damped min-max problem in \eref{LP} is quadratic in both $\u$ and $\Lamm$ variables and its optimum solution is given by the following saddle point system:
\begin{equation} \label{penalty}
\begin{pmatrix}
\bold{A}(\bold{m}) & -\frac{1}{\regp}\bold{I} \\
\bold{P}^T\bold{P} & \bold{A}(\bold{m})^T
\end{pmatrix}
\begin{pmatrix}
\ur^{+} \\
\Lamm^{+}
\end{pmatrix}=
\begin{pmatrix}
\bold{b}^* \\
\bold{P}^T\bold{d}^*
\end{pmatrix}.
\end{equation}
Unlike \eref{redused}, this linear system is more difficult to solve. It can be solved using various standard numerical linear algebra concepts, leading to different yet equivalent iterative procedures.
We can use Schur complement to write its solution explicitly as
\begin{subequations} \label{Explicit_uv}
\begin{align} 
\up^{+}&=\left(\bold{P}^T\bold{P}+\regp\bold{A}^{T}\bold{A}\right)^{\!\!-1}(\bold{P}^T\bold{d}^*+\regp \bold{A}^T\bold{b}^* ), \label{Explicit_u}\\
\Lamm^{+} &= \left(\bold{I} + \frac{1}{\regp} \bold{A}^{-T}\!\bold{P}^T\!\bold{P} \bold{A}^{-1} \right)^{\!\!-1}\!\!\bold{A}^{-T}\bold{P}^T(\bold{d}^* - \bold{P} \bold{A}^{-1}\!\bold{b}^* ) \label{Explicit_v1}\\
& = \bold{A}^{-T}\!\bold{P}^T \left(\bold{I} + \frac{1}{\regp} \bold{P} \bold{A}^{-1}\!\! \bold{A}^{-T}\!\bold{P}^T\right)^{\!\!-1}\!\!(\bold{d}^* - \bold{P} \bold{A}^{-1} \bold{b}^* ). \label{Explicit_v2}
\end{align}
\end{subequations}
where $\bold{A}\equiv\bold{A}(\bold{m})$.
The equality of equations \eqref{Explicit_v1} and \eqref{Explicit_v2} can be found in \citet{Guttman_1946_EMC}. However, while both equations are equivalent, \eref{Explicit_v1} can be more challenging to compute due to the need for the inversion of $(\bold{I} + \frac{1}{\regp} \bold{A}^{-T}\bold{P}^T\bold{P} \bold{A})$, which is of size $N\times N$. On the other hand, \eref{Explicit_v2} involves the inversion of a significantly smaller matrix $(\bold{I} + \frac{1}{\regp} \bold{P} \bold{A}^{-1} \bold{A}^{-T}\bold{P}^T)$ of size $M\times M$.

Due to the complexity of matrix inversion in both cases, it is possible to compute either $\up^{+}$ or $\Lamm^{+}$ explicitly and then use the first equation in \eref{penalty} to compute the other variable. This approach leads to two different strategies for solving the equations: the wavefield-oriented approach and the multiplier-oriented approach.
In the wavefield-oriented approach, we first compute $\up^{+}$ explicitly using \eref{Explicit_u}, and then we utilize the first equation in \eref{penalty} to solve for $\Lamm^{+}$.
Conversely, in the multiplier-oriented approach, we initially compute $\Lamm^{+}$ explicitly using either equation \eqref{Explicit_v1} or \eqref{Explicit_v2}. Subsequently, we employ the first equation in \eref{penalty} to determine $\up^{+}$.

\subsubsection{Wavefield-oriented penalty method}
The wavefield-oriented penalty method results in the following iteration:
\begin{subequations}\label{WMI_penalty} 
\begin{align} 
&\up^{+}=\left(\bold{P}^T\bold{P}+\regp \bold{A}(\bold{m})^T\!\bold{A}(\bold{m})\right)^{-1}\!\!(\bold{P}^T\bold{d}^*+\regp \bold{A}(\bold{m})^T\bold{b}^*),\label{WMI_penaltya}\\
&\Lamm^{+} = \regp (\bold{A}(\bold{m})\up^{+}-\bold{b}^*),\label{WMI_penaltyb}\\
& {m}^+={m}-\alpha \frac{\langle \Lamm^+, \partial_{tt} \u^{+}\rangle_t }{\langle \partial_{tt} \u^{+}, \partial_{tt} \u^{+}\rangle_t} . \label{WMI_penaltyc}
\end{align}
\end{subequations} 

\subsubsection{Multiplier-oriented penalty method}
The multiplier-oriented penalty method results in the following iteration:
\begin{subequations}\label{MWI_penalty}
\begin{align} 
&\Lamm^{+} =  \bold{A}^{-T}\!\bold{P}^T \left(\bold{I} + \frac{1}{\regp} \bold{P} \bold{A}^{-1}\!\! \bold{A}^{-T}\!\bold{P}^T\right)^{-1}\!\!(\bold{d}^* \!- \bold{P} \bold{A}^{-1} \bold{b}^* ),  \label{MWI_penaltya}\\
&\up^{+} = \bold{A}^{-1}(\bold{b}^* + \frac{1}{\regp}\Lamm^{+}),\label{MWI_penaltyb}\\
& {m}^+={m}-\alpha \frac{\langle \Lamm^+, \partial_{tt} \u^{+}\rangle_t }{\langle \partial_{tt} \u^{+}, \partial_{tt} \u^{+}\rangle_t} . \label{MWI_penaltyc}
\end{align}
\end{subequations}
where $\bold{A}\equiv\bold{A}(\bold{m})$.


\subsubsection{Wavefield-oriented vs. multiplier-oriented approaches}

Both wavefield-oriented and multiplier-oriented approaches are two alternative (equivalent) strategies for simultaneously computing the wavefield and multiplier, solving a saddle point system arising from Lagrangian-based FWI. The wavefield-oriented approach emphasizes wavefield reconstruction as its primary focus, while in the multiplier-oriented approach the vector of Lagrange multipliers is computed first. Each approach offers distinct perspectives and trade-offs in terms of computational efficiency and complexity. The choice between the two depends on the specific requirements and constraints of the problem at hand.

In the wavefield-oriented approach, the main challenge lies in inverting the associated Hessian matrix $(\bold{P}^T\bold{P}+\regp \bold{A}^T\bold{A})$ for updating the wavefield. In the context of time domain applications, this matrix can be difficult to invert due to the lack of a suitable time-stepping algorithm \citep{Wang_2016_FIR, Rizzuti_2021_DFW, Gholami_2022_EFW, Lin_2023_FWR}. For small to medium scale problems in frequency domain, direct factorization methods like LU triangular factorization can be used to invert the matrix efficiently when applied to multiple sources. However, for large-scale applications, iterative methods like CG become more feasible, although appropriate preconditioning is necessary to handle the severe ill-conditioning of the Hessian matrix.

On the other hand, the multiplier-oriented approach offers the advantage of computing the wavefield using standard time-stepping methods, as long as the multiplier is available. In this approach, the computational bottleneck lies in the computation of the LS multipliers, which involves inverting the dense Hessian matrix $(\bold{I} + \frac{1}{\regp} \bold{P} \bold{A}^{-1} \bold{A}^{-T}\bold{P}^T)$. However, the size of this matrix is the same as the data size \citep[see Appendix A in][]{Gholami_2022_EFW}. While the explicit construction of this matrix requires solving $N_r$ adjoint wave equations, its direct inversion is possible only for frequency domain applications. In time domain, approximate solutions can be obtained using iterative solvers such as preconditioned CG.
In this case, it is not necessary to explicitly build the Hessian matrix but only to compute the product of it by a vector, which incurs the cost of two additional wave-equation solves per CG iteration.

\subsubsection{Connection with WRI}
The iteration \eqref{WMI_penalty} simplifies into the WRI method proposed by \citet{VanLeeuwen_2013_MLM} when the step length $\alpha$ is set as the inverse of the penalty parameter, i.e., $\alpha=1/\regp$. By substituting the expression for the multiplier $\Lamm^{+}$ from \eref{WMI_penaltyb} into the model update formula in \eref{WMI_penaltyc}, we derive a two-step WRI iteration as shown in \eref{WRI}.
However, it is important to note that WRI inherently does not satisfy the wave equation constraint. Therefore, the penalty parameter $\regp$ needs to be progressively increased during iterations to approach convergence, typically by employing a parameter continuation strategy. Yet, this strategy leads to an exacerbated ill-conditioning issue, which is a fundamental limitation of the WRI method. This arises from the fact that a larger $\regp$ results in a smaller step length, causing the algorithm to converge slowly.

This issue can be more effectively demonstrated using the formulation \eqref{MWI_penalty}, which is equivalent to \eqref{WMI_penalty}. By choosing $\alpha=1/\regp$ and introducing a change of variable $\Lam=\frac{1}{\regp}\Lamm$ in \eref{MWI_penalty}, we arrive at the following multiplier-oriented variant of WRI:
 \begin{subequations}\label{MWI_penalty2}
\begin{align} 
&\Lam^{+} =  \bold{A}^{-T}\!\bold{P}^T \left(\regp\bold{I} + \bold{P} \bold{A}^{-1}\!\! \bold{A}^{-T}\!\bold{P}^T\right)^{-1}\!\!(\bold{d}^* \!- \bold{P} \bold{A}^{-1} \bold{b}^* ),  \label{MWI_penalty2a}\\
&\up^{+} = \bold{A}^{-1}(\bold{b}^* + \Lam^{+}),\label{MWI_penalty2b}\\
& {m}^+={m}- \frac{\langle \Lam^+, \partial_{tt} \u^{+}\rangle_t }{\langle \partial_{tt} \u^{+}, \partial_{tt} \u^{+}\rangle_t}, \label{MWI_penalty2c}
\end{align}
\end{subequations}
This formulation demonstrates that the parameter $\regp$ effectively controls the amount of damping applied to the (scaled) multiplier $\Lam^{+}$. As the value of $\regp$ increases, the damping effect intensifies, potentially damping the multiplier to zero. This excessive damping causes ill-conditioning, which reflects in a slow convergence rate.

One way to address this challenge is by relaxing the constraint $\alpha=1/\regp$, allowing for greater flexibility in controlling the algorithm's convergence. Moreover, a robust solution to this issue is offered by the AL formulation. The AL method operates effectively with moderate values of the penalty parameter, avoiding the need for excessively large values. This preserves algorithm well-posedness and enhances convergence behavior.
Maintaining algorithm well-posedness comes at the cost of retaining and updating estimates of the Lagrange multipliers in each iteration. While this introduces increased memory requirements, practical strategies can be employed to efficiently manage this memory demand \citep{Gholami_2022_EFW}.

\subsection{Regularization by the AL method}
An excellent method for regularization of the Lagrangian problem is using the proximal regularization, that allows prior information of the multiplier to be introduced into the objective function \citep{Rockafellar_1976_MOA,Gill_2012_PDA,Wright_2022_ODA}:
\begin{equation}  \label{LAL}
(\u^{+},\Lamm^{+})=\arg \minimax_{\u,\Lamm} ~\mathcal{L}(\bold{m},\u,\Lamm) - \frac{1}{2\regp}\|\Lamm-\LamAL\|_2^2.
\end{equation} 
where $\LamAL$ is a prior multiplier estimate. 
The physical interpretation of this strategy is simple: when maximizing the Lagrangian function, one wishes a multiplier $\Lamm$ that is not too far from the a priori multiplier $\LamAL$. The constraints that the updated multiplier cannot be too far from some a priori estimate is necessary for avoiding ill-posedness and building a stable iteration \citep{Tarantola_1986_SNL}. In the proximal regularization, the prior estimate $\LamAL$ varies with iteration \citep{Rockafellar_1976_MOA}. 
It may be selected as the best estimation of the original multiplier at the previous iteration, 
\begin{equation} \label{AL_das}
\LamAL^{+}=\LamAL+\regp(\bold{A}(\bold{m}^{+})\um^{+}-\bold{b}^*). 
\end{equation}
The proximal Lagrangian approach is closely related to the AL method, but it provides a more intuitive and easier-to-interpret formulation. 
Note that the objective function in \eref{LAL} is quadratic in $\Lamm$ and has a simple solution $\Lamm =\LamAL+ \regp (\bold{A}(\bold{m})\ur - \bold{b}^*)$. Substituting this expression into the objective function eliminates the primal Lagrange multiplier from the optimization variables, giving the standard form of the AL function in \eref{AL_obj}.

The optimum solution of the objective function in \eref{LAL} is given by the following saddle point system:
\begin{equation} \label{AL}
\begin{pmatrix}
\bold{A}(\bold{m}) & -\frac{1}{\regp}\bold{I} \\
\bold{P}^T\bold{P} & \bold{A}(\bold{m})^T
\end{pmatrix}
\begin{pmatrix}
\ur^{+} \\
\Lamm^{+}
\end{pmatrix}=
\begin{pmatrix}
\bold{b}^*-\frac{1}{\regp}\LamAL \\
\bold{P}^T\bold{d}^*
\end{pmatrix}.
\end{equation}
This system is essentially the same as that in \eref{penalty}, with the only difference being the substitution of the physical source $\bold{b}^*$ at the right-hand side in the \eref{penalty} with the extended source $\bold{b}^*-\frac{1}{\regp}\LamAL$ in \eref{AL}.
Therefore, we can employ the two approaches presented above to solve it.


\subsubsection{Wavefield-oriented AL method}
Applying this approach gives us the following iteration:
\begin{subequations}\label{IRWRI}
\begin{align} 
&\up^{+} =\left(\bold{P}^T\bold{P}+\regp \bold{A}(\bold{m})^T\bold{A}(\bold{m})\right)^{-1}(\bold{P}^T\bold{d}^*+\regp \bold{A}(\bold{m})^T\bold{b}^*-\bold{A}(\bold{m})^T\LamAL),\label{IRWRIa}\\
&\Lamm^{+} = \LamAL +\regp (\bold{A}(\bold{m})\up^{+}-\bold{b}^*),\label{IRWRIb}\\
& {m}^+={m}-\alpha \frac{\langle \Lamm^+, \partial_{tt} \u^{+}\rangle_t }{\langle \partial_{tt} \u^{+}, \partial_{tt} \u^{+}\rangle_t}, \label{IRWRIc}\\
& \LamAL^{+}=\LamAL+\regp(\bold{A}(\bold{m}^{+})\um^{+}-\bold{b}^*). \label{IRWRId}
\end{align}
\end{subequations} 
If we employ the step length $\alpha=1/\regp$ in \eref{IRWRIc} to update the model, then this iteration will be equivalent to IR-WRI \citep{Aghamiry_2019_IWR}. IR-WRI faces a similar challenge as WRI regarding the wavefield reconstruction step, \eref{IRWRIa}.

\subsubsection{Multiplier-oriented AL method}
The multiplier-oriented approach solves the AL problem by the following iteration:
\begin{subequations}\label{MOALM}
\begin{align} 
&\Lamm^{+} =  \bold{A}(\bold{m})^{-T}\!\bold{P}^T \left(\bold{I} + \frac{1}{\regp} \bold{P} \bold{A}(\bold{m})^{-1} \bold{A}(\bold{m})^{-T}\bold{P}^T\right)^{-1}(\bold{d}^* - \bold{P} \bold{A}(\bold{m})^{-1} \bold{b}^* + \frac{1}{\regp} \bold{P} \bold{A}(\bold{m})^{-1}\LamAL), \label{MOALMa}\\
&\up^{+} = \bold{A}(\bold{m})^{-1}(\bold{b}^* + \frac{1}{\regp}\Lamm^{+}-\frac{1}{\regp}\LamAL),\label{MOALMb} \\
& {m}^+={m}-\alpha \frac{\langle \Lamm^+, \partial_{tt} \u^{+}\rangle_t }{\langle \partial_{tt} \u^{+}, \partial_{tt} \u^{+}\rangle_t}, \label{MOALMc}\\
& \LamAL^{+}=\LamAL+\regp(\bold{A}(\bold{m}^{+})\um^{+}-\bold{b}^*). \label{MOALMd}
\end{align}
\end{subequations}

\subsubsection{Penalty vs. AL methods}
In comparison to the penalty-based methods the AL-based iterations involve an additional secondary Lagrange multiplier $\LamAL$. This multiplier plays a crucial role in enhancing the convergence of the algorithm.
Specifically, in AL-based iterations, the wave equation depicted in \eref{MOALMb} incorporates three distinct types of sources. The known physical source, $\bold{b}^*$, serves as the primary source, generating the background wavefield. The Lagrange multiplier $\Lamm$ acts as a secondary volume source, contributing to the low-order scattering part of the wavefield. The Lagrange multiplier $\LamAL$ acts as a tertiary volume source, correcting the defects in $\Lamm$ and contributing to the high-order scattering part of the wavefield. This tertiary volume source, with its significant role in conditioning the iteration, is absent in penalty-based iterations.

The penalty-based iterations have been extensively analyzed in the framework of WRI by \cite{Operto_2023_FWI}, but the numerical results were produced using the AL iteration, including this tertiary volume source term. Although there are similarities and relationships between the two approaches, they are not the same and can yield completely different results. 
To better comprehend the similarities and differences between these methods and the mechanisms employed by each of them in solving FWI, we need to establish the physical meaning and precisely analyze the actions of each term in the optimization process.
In the subsequent section, we will provide a deeper understanding of the physical meaning of both $\Lamm$ and $\LamAL$ and their roles in improving the inversion results within a general framework that reveals the interconnections between different presented methods.

\begin{figure*}
\center
\includegraphics[angle=0,scale=0.9]{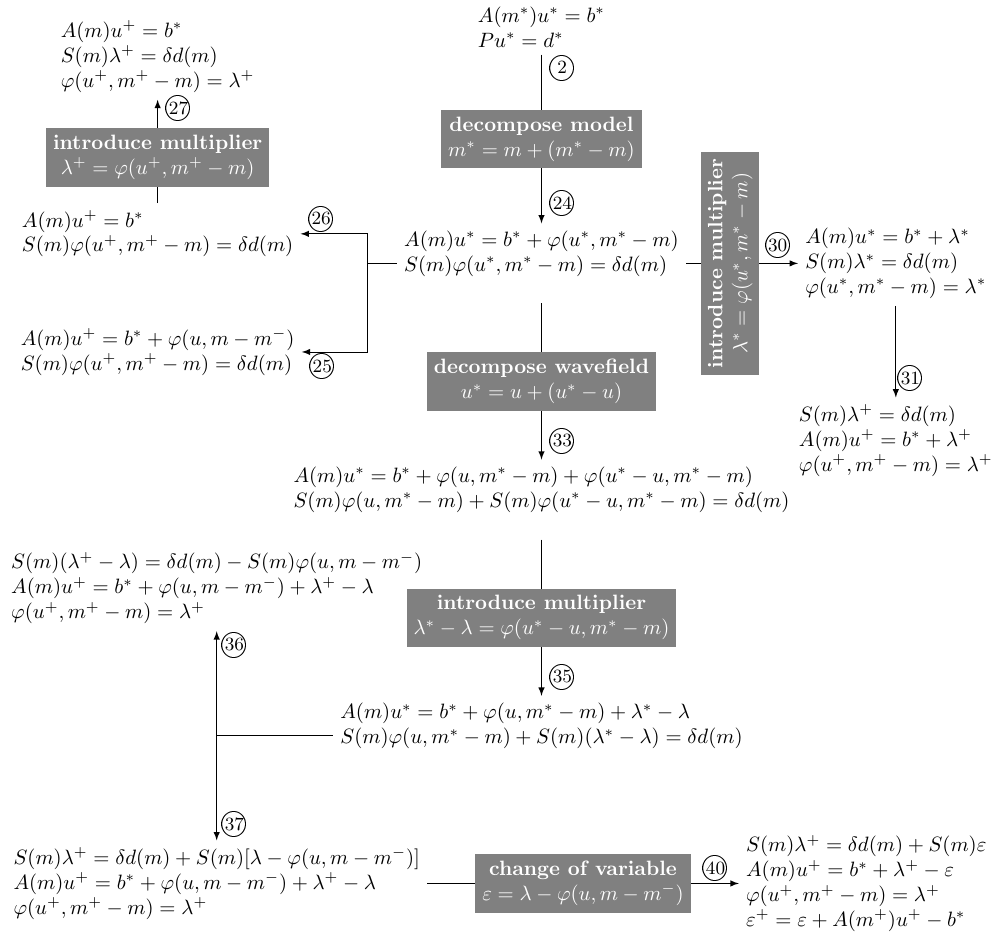}
\caption{Interconnections between iterative methods for solving FWI equations. The circled numbers represent the corresponding equation numbers in the text. Each iteration consists of solving a series of equations sequentially, updating the variables denoted with $\bullet^+$ using the background values $\bullet$ and the previous values $\bullet^-$.
}
\label{chart}
\end{figure*}

\section{The basic iterations, mechanisms, and physical interpretations} \label{physical_interp}
In this section, the Lagrangian based algorithms discussed in the previous section are further developed and examined using basic iterative linearization methods and scattering concepts. Using this formulation, we gain a deeper physical understanding of the role of Lagrange multipliers in the reconstruction process. We can derive the same algorithms obtained through Lagrangian formulation but from a scattering point of view and simple alternating iterative methods. This offers a unique perspective that connects the optimization theory with the physics of wave scattering. 
We present the iterations in a progressive manner, starting from their simplest form and gradually advancing towards more sophisticated and refined versions. Throughout the development, we discuss the advantages and disadvantages of each iteration, highlighting their evolution and improvements. Importantly, we establish clear connections between the generated iterations and their corresponding counterparts in Lagrangian optimization, whenever necessary. By taking this comprehensive approach, we enable a thorough comparison and evaluation of the various Lagrangian based methods discussed earlier. This analysis allows us to assess their performance, understand their underlying mechanisms, and identify their strengths and limitations from the scattering theory perspective.


Let us decompose the true model $\bold{m}^*$ into the known background model $\bold{m}$ and the unknown perturbation $\bold{m}^*-\bold{m}$, then we can split the nonlinear operator $\bold{A}(\bold{m}^*)\bold{u}^*$ into linear and bilinear terms, $\bold{A}(\bold{m})\bold{u}^*$ and $\scs(\bold{u}^*,\bold{m}^*-\bold{m})$, respectively, as 
\begin{equation} 
\bold{A}(\bold{m}^*)\bold{u}^* = \bold{A}(\bold{m})\bold{u}^*-\scs(\bold{u}^*,\bold{m}^*-\bold{m}),
\end{equation}
\citep{Schuster_1985_HBB}. Accordingly, \eref{waveEq} can be written in the form of the Lippmann–Schwinger integral equation \citep{Lippmann_1950_VPS} 
\begin{equation} \label{Lippmann}
\bold{A}(\bold{m})\bold{u}^* =\bold{b}^* +\scs(\bold{u}^*,\bold{m}^*-\bold{m}).
\end{equation}
The bilinear function $\scs$, defined as
\begin{equation} \label{scs}
\scs(\bold{u}^*,\bold{m}^*-\bold{m}) =\left(-\frac{\partial \bold{A(m)}}{\partial \bold{m}}\u^*\right) (\bold{m}^*-\bold{m}),
\end{equation}
is the scattering source, also known as ``secondary Born source" \citep{Tarantola_1988_TBI} or ``contrast source" \citep{vandenBerg_1997_CSI}. It is created by product of the perturbation, $\bold{m}^*-\bold{m}$, and the wavefield, and as a result, it may extend the entire space-time domain.
When both arguments are considered simultaneously, $\scs$ is nonlinear, but when only one argument is considered, it is linear.

The following data equation connects the scattered data or data residual $\delta\bold{d(m)}=\bold{d}^*-\bold{P}\bold{A(m})^{-1}\bold{b}^*$ to the scattering source linearly by combining \eref{dataEq} and \eref{Lippmann}: 
\begin{equation} \label{delta_d}
\bold{S}(\bold{m})\scs(\bold{u}^*,\bold{m}^*-\bold{m})=\delta\bold{d(m)}.
\end{equation}
where $\bold{S(m)}=\bold{P}\bold{A(m})^{-1}$.
This equation represents a forward modeling operator that calculates the scattered data resulting from a given scattering source $\scs$. However, it is still nonlinear due to the dependence of the incident wavefield $\bold{u}^*$ on the model parameters. We can linearize the equation around the background model $\bold{m}$ by treating the wavefield as an independent variable. The accuracy of this linearization is influenced by the accuracy of the wavefield, and an exact linearization can be achieved if the exact wavefield is known.

We use the splitting above to yield a new pair of coupled equations, as given in \eref{Pair2}, which are alternatives to the equations in \eref{Pair1}.
\begin{subequations} \label{Pair2}
\begin{align}
&\bold{A}(\bold{m})\bold{u}^* =\bold{b}^* +\scs(\bold{u}^*,\bold{m}^*-\bold{m}),\label{Pair2a}\\
&\bold{S}(\bold{m})\scs(\bold{u}^*,\bold{m}^*-\bold{m})=\delta\bold{d}(m). \label{Pair2b}
\end{align}
\end{subequations}
It should be noted that the two sets of equations in \eqref{Pair1} and \eqref{Pair2} are equivalent. This means that the optimal solution $(\bold{m}^*,\bold{u}^*)$  obtained from either set of equations also optimally solves the other set.
The coupled equations in \eqref{Pair2} serve as the foundation to develop a variety of alternating iterations. 
As we will show in what follows, a number of interesting iterations, including the above Lagrangian-based algorithms, can be derived from these equations (see Figure \ref{chart}). The number and type of chosen variables, as well as the linear equations between them, differ. The generated iterations may be equivalent in that they arrive at the same solution. However, the solution process may be significantly affected by the chosen formulation.

\subsection{The Gauss-Seidel iteration}
The coupled equations presented in \eref{Pair2} pose a challenge for simultaneous solution due to their inherent nonlinearity. To overcome this challenge, an iterative decoupling approach can be employed. This method breaks down the problem into subproblems, with each subproblem focusing on updating only one variable by solving a single equation using information from the previous iteration. 

Let us start with an initial point $(\bold{m}_0,\bold{u}_0)$, set $\bold{m}_{-1}=\bold{m}_0$, and solve \eref{Pair2a} for the new wavefield $\bold{u}_1$ using $\scs(\bold{u}_0,\bold{m}_0-\bold{m}_{-1})$ as an estimate of the scattering source. We determine a new model $\bold{m}_1$ such that $\scs(\bold{u}_1,\bold{m}_1-\bold{m}_0)$ solves \eref{Pair2b} after obtaining the new wavefield $\bold{u}_1$.
Subsequently, with $(\bold{m}_1,\bold{u}_1)$, the process is repeated to obtain $(\bold{m}_2,\bold{u}_2)$.
 This Gauss-Seidel procedure is repeated until convergence is achieved. The general form of the iteration is as follows.
\begin{subequations} \label{GS}
\begin{align}
&\bold{A}(\bold{m})\bold{u}^+ =\bold{b}^* +\scs(\bold{u},\bold{m}-\bold{m}^-),\label{GSa}\\
&\bold{S}(\bold{m})\scs(\bold{u}^+,\bold{m}^+-\bold{m})=\delta\bold{d}(m). \label{GSb}
\end{align}
\end{subequations}
The concept behind this iteration is straightforward. By utilizing an approximate wavefield and model perturbation, we can construct an approximate scattering source in the Lippmann–Schwinger equation. A natural choice is to use the most recent estimates, $\scs(\bold{u},\bold{m}-\bold{m}^-)$, which leads to an approximation of the wavefield  $\bold{u}^+$ using \eref{GSa} and then model $\bold{m}^+$ using \eref{GSb}. 

In the Gauss-Seidel iteration, the nonlinearity of the problem is effectively addressed through two key factors. Firstly, the model is updated at each iteration, allowing for the gradual inclusion of scattering terms into the model parameters. Secondly, it accounts for remaining multiscattering effects encoded in $\scs(\bold{u},\bold{m}-\bold{m}^{-})$, which enables the iteration to capture higher-order scattering contributions, akin to a Born series solution. By employing the previous estimates of the wavefield and perturbation, this iterative approach represents a higher-order approximation, resulting in an enhanced estimation of the wavefield and subsequently, the model parameters. This capability allows for better convergence and improved inversion results.

However, there are two significant challenges to implementing the Gauss-Seidel iteration effectively. The first challenge pertains to computational costs, as solving \eref{GSb} for the updated model is very challenging due to the large size of the associated kernel matrix which is source dependent. 
The second challenge lies in memory requirements. Storing and accessing the previous wavefield estimates can lead to substantial memory usage, which may hinder the applicability of the method to large-scale problems. 

\subsection{The Gauss-Newton iteration}
Even though the Gauss-Newton based FWI is particularly recognized as a gradient based iteration for minimization of the reduced objective function in  \eref{FWI_obj} \citep{Pratt_1998_GNF, Metivier_2017_TRU}, the resulting iteration can simply be obtained from the Gauss-Seidel iteration, \eref{GS}, by omitting the scattering source term $\scs$ from the right-hand side of \eref{GSa}: 
\begin{subequations} \label{GN}
\begin{align}
&\bold{A}(\bold{m})\bold{u}^+ =\bold{b}^*,\label{GNa}\\
&\bold{S}(\bold{m})\scs(\bold{u}^+,\bold{m}^+-\bold{m})=\delta\bold{d}(m). \label{GNb}
\end{align}
\end{subequations}
The wavefield is first computed by solving the wave equation, \eref{GNa}, with the physical source, starting from the background model $\bold{m}$.
We then solve the data equation \eqref{GNb} to update the model using the computed wavefield. 
Gauss-Newton method addresses the nonlinearity of the problem only through iteration by updating the background model. 
The trick used in this iteration substantially reduces the memory requirements associated with storage of the wavefields; however, it suffers from approximations in the \eref{Pair2a} of replacing the scattering source with zero. In the case of a strong multiple scattering effect, such an approximation is inaccurate which may cause instability. However, when the background model is sufficiently accurate, this term is small and this approximation can be acceptable. 

\subsection{The Split Gauss-Newton iteration}
The Gauss-Newton method encounters a significant computational challenge when solving \eref{GNb}. Directly solving this linear system can be computationally demanding and may not be practical for large-scale problems. To overcome this issue, several strategies can be employed. One option is to use an iterative method, such as the CG method, to approximate the solution of \eref{GNb} \citep{Metivier_2017_TRU}.
Another approach to address the computational challenge is to take advantage of the split structure of the kernel matrix and decompose the solution procedure into two sequential steps,  using the same strategy as solving a linear system of equations by LU decomposition  \citep{Gholami_2022_OCN}. In this alternative approach, we introduce an auxiliary variable (multiplier) $\Lam^+=\scs(\bold{u}^+,\bold{m}^+-\bold{m})$, and then we split \eref{GNb} into the multiplier equation, \eref{SGNb}, and constitutive equation,  \eref{SGNc}, which are solved in sequence, leading to the split Gauss-Newton iteration:
\begin{subequations} \label{SGN}
\begin{align}
&\bold{A}(\bold{m})\bold{u}^+ =\bold{b}^*,\label{SGNa}\\
&\bold{S}(\bold{m})\Lam^+=\delta\bold{d}(m), \label{SGNb}\\
&\scs(\bold{u}^+,\bold{m}^+-\bold{m}) = \Lam^+. \label{SGNc}
\end{align}
\end{subequations}
Using the LS approach to solve the  constitutive equation gives
\begin{equation} \label{m_update}
\minimize_{\bold{m}^+}  \|\scs(\u^{+},\bold{m}^+-\bold{m})-\Lam^{+}\|_2^2 \quad \Longrightarrow \quad {m}^+={m}- \frac{\langle \Lam^+, \partial_{tt} \u^{+}\rangle_t }{\langle \partial_{tt} \u^{+}, \partial_{tt} \u^{+}\rangle_t}.
\end{equation}
We can see that the model update in this case is very similar to that in the Lagrangian based methods, \eref{fwi_dm}. The step length here is unity and $\Lamm^{+}$ is replaced with $\Lam^{+}$.
Interestingly, when we examine the iteration \eref{SGN}, we can notice that it can be reduced to the Lagrangian-based iteration in \eref{redused2} if we estimate the multipliers by using an adjoint operator to approximately solve the multiplier equation, \eref{SGNb}. The step length $\regp$ in \eref{redused2} partially accounts for these approximations. In this context, the Lagrange multipliers $\Lam^{+}$ and $\Lamm^{+}$ are related by the equation $\Lam^{+}=(1/\regp)\Lamm^{+}$.



\subsection{The Split Gauss-Seidel iteration}
By introducing the auxiliary variable 
\begin{equation}
\Lam^*=\scs(\bold{u}^*,\bold{m}^*-\bold{m}),
\end{equation}
we can split the equations in \eref{Pair2} as
\begin{subequations} \label{split}
\begin{align}
&\bold{A}(\bold{m})\bold{u}^* =\bold{b}^* +\Lam^*,\label{splita}\\
&\bold{S}(\bold{m})\Lam^*=\delta\bold{d}(m), \label{splitb}\\
&\scs(\bold{u}^*,\bold{m}^*-\bold{m})=\Lam^*. \label{couple}
\end{align}
\end{subequations}
If we reverse the order of the first and second equations and perform an alternating solve, we get \eref{SGS}, a split form of the Gauss-Seidel iteration in \eref{GS}.
\begin{subequations} \label{SGS}
\begin{align}
&\bold{S(m)}\Lam^{+}=\delta\bold{d}(m), \label{SGSb}\\
&\bold{A}(\bold{m})\bold{u}^{+} =\bold{b}^* +\Lam^{+},\label{SGSa}\\
&\scs(\u^{+},\bold{m}^{+}-\bold{m})=\Lam^{+}.\label{SGSc}
\end{align}
\end{subequations}
In this iteration, the optimum $\Lam^+$ is the true scattering source, which is an optimum solution of \eref{SGSb}. Thus, if matrix $\bold{S}$ is invertible, we can compute $\Lam^*$ precisely from this equation in Step 1. Then, in step 2, the true wavefields and at step 3 converge to the true model in a single iteration. The main difficulty is that $\bold{S}$ is not invertible, which demands a more complicated iterative process to solve the problem.
First, the multiplier equation \eqref{SGSb} is solved for the updated scattering source $\Lam^{+}$ and then the wavefield is obtained by solving \eref{SGSa} using $\Lam^{+}$. Finally, model $\bold{m}^{+}$ is recovered by factoring the scattering source estimate into the product of the wavefield and model perturbation. 
Since $\bold{A}$ is invertible we can use the expression of the multiplier from \eref{SGSa} in \eref{SGSc} to get
$\scs(\u^{+},\bold{m}^{+}-\bold{m})=\bold{A}(\bold{m})\bold{u}^{+} -\bold{b}^*$. Therefore, in essence, the split Gauss-Seidel iteration, \eref{SGS}, updates the model parameters using the Lippmann–Schwinger equation unlike the Gauss-Seidel iteration in \eref{GS} which uses the data equation for this update.

This split iteration may have the following advantages: (1) it does not require storage of the previous wavefields. (2) It can be more efficient for inversion of multiple source data using the same spread acquisition, because we only need to construct and invert the dense matrix $\bold{S}$ once and reuse it for various sources because it is source-independent. (3) The bilinear operator $\csc$ and the matrix $bold{S}$ have different conditioning; the former is rank defficient and may need strong regularization, whilst the latter is well posed thus inversion in the split form permits more flexibility in the suitable regularization performed at each level. 


%
The split Gauss-Seidel iteration can be equivalent to the multiplier-oriented penalty method in  \eref{MWI_penalty}. This equivalence arises for the case $\alpha=1/\regp$ and when we solve the multiplier equation in \eref{SGSb} using a damped LS approach with a damping parameter $\regp$. In this context, the associated Lagrange multipliers are related by the equation $\Lam^{+}=(1/\regp)\Lamm^{+}$. This equivalence provides evidence that the quadratic penalty methods like WRI utilizes a first-order Lagrange multiplier, which represents the scattering source that best fits the scattered data in a LS sense.
%
%

%

\subsection{Approximations and errors} \label{erros}
The iterations described by eqs. \eqref{GS}, \eqref{GN}, \eqref{SGN} and \eqref{SGS}  exhibit several sources of error that are not adequately considered. It is important to recognize and address these errors to ensure accuracy and reliability. Some notable errors include the following.
\begin{enumerate}
\item Approximation error of the estimated scattering sources. The estimation of the scattering source involves certain approximations and assumptions. It is crucial to quantify and minimize these errors to obtain more reliable results. 

\item Regularization errors. Regularization techniques are employed to stabilize the inversion process and handle ill-posed equations during updating of the model and multipliers. However, regularization introduces additional errors due to the choice of regularization parameters and assumptions regarding the statistics of the solution. 

\item Simplifications in solving the equations: The iterative process often involves simplifications and assumptions in solving the system of equations involved in updating the model and multipliers. These simplifications may neglect certain terms or assume certain conditions that can lead to inaccuracies in the solutions. 

\end{enumerate}
In the following section, we describe refined iterations that appropriately account for these errors. 

\subsection{Accounting for high order scattering terms and errors}

Let us decompose the true wavefield $\bold{u}^*$ into the incident wavefield $\bold{u}$ plus the scattered wavefield $\bold{u}^*-\bold{u}$.
It should be noted that $\bold{u}$ is not necessarily the background wavefield generated by the physical source in the background model. Indeed, this is a low-order approximation of the true wavefield $\bold{u}^*$. Thus, it may satisfy the wave equation in $\bold{m}$ only approximately, that is, $\bold{A}(\bold{m})\bold{u}\approx \bold{b}^*$. If $\bold{u}$ is the background wavefield satisfying $\bold{A}(\bold{m})\bold{u}= \bold{b}^*$, then $\bold{u}^*-\bold{u}$ is the total scattered wavefield; otherwise, it captures only the high-order part of the scattered field.

This decomposition allows us to separate the dominant (low-order) scattering term from the residual (high-order) terms as follows:
\begin{equation} \label{scatterers}
\scs(\bold{u}^*,\bold{m}^*-\bold{m})=\scs(\bold{u},\bold{m}^*-\bold{m})+\scs(\bold{u}^*-\bold{u},\bold{m}^*-\bold{m}),
\end{equation}
where $\scs(\bold{u},\bold{m}^*-\bold{m})$ represents the dominant part of the scattering source and $\scs(\bold{u}^*-\bold{u},\bold{m}^*-\bold{m})$ represents the change in $\scs$ due to the error in the wavefield. 
Using this equation in \eref{Pair2}, we create another set of equations to be solved by alternating iterations:
\begin{subequations} \label{Pair3}
\begin{align}
&\bold{A}(\bold{m})\bold{u}^* =\bold{b}^* +\scs(\bold{u},\bold{m}^*\!\!-\bold{m})+\scs(\bold{u}^*\!\!-\bold{u},\bold{m}^*\!\!-\bold{m}), \label{Pair3a}\\
&\bold{S}(\bold{m})\scs(\bold{u},\bold{m}^*\!\!-\bold{m})+\bold{S}(\bold{m})\scs(\bold{u}^*\!\!-\bold{u},\bold{m}^*\!\!-\bold{m})=\delta\bold{d}(\bold{m}). \label{Pair3b}
\end{align}
\end{subequations}
In this form, the data residuals are decomposed into a low-order scattered data component (the first term in \eref{Pair3b}) and the remaining multiscattered part (the second term in \eref{Pair3b}).  

Decomposition of the wavefield also allows us to define a perturbation for the Lagrange multiplier:
\begin{equation} \label{dLambda}
\Lam^*-\Lam=\scs(\bold{u}^*-\bold{u},\bold{m}^*-\bold{m}).
\end{equation}
Utilizing \eref{dLambda} in \eref{Pair3}, we obtain the following equations:
\begin{subequations} \label{Pair33}
\begin{align}
&\bold{A}(\bold{m})\bold{u}^* =\bold{b}^* +\scs(\bold{u},\bold{m}^*-\bold{m})+\Lam^*-\Lam, \label{Pair33a}\\
&\bold{S}(\bold{m})\scs(\bold{u},\bold{m}^*-\bold{m})+\bold{S}(\bold{m})(\Lam^*-\Lam)=\delta\bold{d}(\bold{m}). \label{Pair33b}
\end{align}
\end{subequations}
To solve these equations using alternating iterations, the following steps are followed, requiring $\Lam, \bold{u}$, and $\bold{m}^-$:
\begin{enumerate}
    \item Approximate the model perturbation $\bold{m}^*-\bold{m}$ in the scattering source with the previous value, $\bold{m}-\bold{m}^{-}$.
    \item Estimate the multiplier $\Lam^+$ by solving \eref{Pair33b}.
    \item Use these estimates to construct the scattering source and solve \eref{Pair33a} to estimate the incident wavefield $\bold{u}^+$.
    \item Finally, generate the updated model by solving the constitutive equation using $\Lam^+$ and $\bold{u}^+$.
\end{enumerate}
The iterative process can be summarized by the following equations:
\begin{subequations} \label{Pair4_iter}
\begin{align}
&\bold{S(m)}(\Lam^{+}-\Lam)=\delta\bold{d}(m)-\bold{S}(m)\scs(\bold{u},\bold{m}-\bold{m}^{-}),\label{Pair4_itera}\\
&\bold{A}(\bold{m})\bold{u}^{+} =\bold{b} +\scs(\bold{u},\bold{m}-\bold{m}^{-})+\Lam^{+}-\Lam_,\label{Pair4_iterb}\\
&\scs(\u^{+},\bold{m}^{+}-\bold{m})=\Lam^{+}. \label{Pair4_iterc}
\end{align}
\end{subequations}
 It begins by subtracting the scattered data $\bold{S(m)}\scs(\bold{u},\bold{m}-\bold{m}^{-})$, computed using the most recent values of $\bold{m}^{-}$, $\bold{m}$, and $\bold{u}$, from the data residual. This subtraction isolates the residual portion that the algorithm was unable to predict in the previous iteration, as shown in \eqref{Pair4_itera}. The differential term $\Lam^{+}-\Lam$ is then used to best fit this remaining multiscattered portion of the data, effectively approximating the corresponding source.
Therefore, this iteration can be viewed as estimating the residual multiscattering source by utilizing the multiplier equation and subsequently adding it to the recently estimated scattering source. It can be considered an improved version of iteration \eqref{GS}, as it considers the errors present in that iteration.
Furthermore, unlike iteration \eqref{SGS}, where the total scattering source is estimated as a whole in each iteration by solving $\bold{S(m)}\Lam^{+}= \delta\bold{d(m)}$, this method primarily focuses on estimating the scattering residuals. This is possible because we already have a reliable estimation of the scattering source from the previous iteration, which allows us to refine and improve the solution iteratively.

Additionally, we can rearrange the terms in \eref{Pair4_iter} to obtain a new interpretation as a refined version of the iteration \eqref{SGS}, as shown in \eref{Pair44_iter}.
\begin{subequations} \label{Pair44_iter}
\begin{align}
&\bold{S(m)}\Lam^{+}=\delta\bold{d}(m)+\bold{S}(m)[\Lam-\scs(\bold{u},\bold{m}-\bold{m}^{-})],\label{Pair44_itera} \\
&\bold{A}(\bold{m})\bold{u}^{+} =\bold{b} +\Lam^{+}+\scs(\bold{u},\bold{m}-\bold{m}^{-})-\Lam_,\label{Pair44_iterb} \\
&\scs(\u^{+},\bold{m}^{+}-\bold{m})=\Lam^{+}. \label{Pair44_iterc} 
\end{align}
\end{subequations}
In this form, the term $\Lam-\scs(\bold{u},\bold{m}-\bold{m}^{-})$ captures the proximity in approximating the scattering source in the previous iteration. There may be residual errors in splitting the multiplier into the model perturbation and the wavefield. These errors in the scattering source estimate are then added back to the data residual before estimating a new multiplier.

Thus, this refined iteration combines elements from the Gauss-Seidel iteration and its split variant, in eqs. \eqref{GS} and \eqref{SGS}, carefully addressing the encountered errors and focusing on refining the scattering sources. The main advantage of this refined iteration is that upon convergence, where $\|\Lam^{+}-\Lam\|= 0$, we obtain the relationship $\bold{S}(m)\scs(\bold{u},\bold{m}-\bold{m}^{-})=\delta\bold{d}(m)$. This remarkable result allows us to approximately solve eqs. \eqref{Pair44_itera} and \eqref{Pair44_iterc} by utilizing regularized formula, while still achieving convergence towards the solution of the inverse problem.

However, it is important to note that the iterations described by eqs. \eqref{Pair4_iter} and \eqref{Pair44_iter} are primarily presented for interpretation purposes and may not be the most efficient for implementation in practice. One limitation is the requirement to store the previous values of both $\bold{u}$ and $\Lam$ in memory. This can be problematic, especially for large-scale problems where memory usage is a concern.
This issue can be partially addressed through a variable change, which leads to a different interpretation of the iteration. Let us introduce the variable $\LamALn$ as the difference between the Lagrange multiplier $\Lam$ and scattering term $\scs(\bold{u},\bold{m}-\bold{m}^-)$,
\begin{equation}
\LamALn=\Lam - \scs(\bold{u},\bold{m}-\bold{m}^-).
\end{equation} 
In fact, $\LamALn$ represents the error in the estimated scattering source term, including the high-order scattering part that was not predicted as well as other errors in solving the equations (see Section \ref{erros}). By using equation \eref{Pair44_iterb} we can estimate the error term $\LamALn^{+}$ in the updated variables:
\begin{equation}
\LamALn^{+}
=\Lam^{+} - \scs(\bold{u}^{+},\bold{m}^{+}-\bold{m})
=\LamALn + \bold{A}(\bold{m}^{+})\bold{u}^{+} - \bold{b}^*.
\end{equation} 
Accordingly, we obtain an equivalent four-step iteration in terms of the error variable:
\begin{subequations} \label{Pair5_iter}
\begin{align}
&\bold{S}(\bold{m})\Lam^{+}=\delta\bold{d(m)}+\bold{S}(\bold{m})\LamALn,\label{Pair5_itera}\\
&\bold{A}(\bold{m})\bold{u}^{+} =\bold{b}^* +\Lam^{+}-\LamALn,\label{Pair5_iterb}\\
&\scs(\u^{+},\bold{m}^{+}-\bold{m})=\Lam^{+}, \label{Pair5_iterc}\\
&\LamALn^{+}=\LamALn + \bold{A}(\bold{m}^{+})\bold{u}^{+}-\bold{b}^*.
\end{align}
\end{subequations}
In this iteration, the only variable that must be stored in memory is the scattering source error $\LamALn$. This modification reduces the memory requirements by 50\% compared with previous iteration.
This four-step iteration can be interpreted as a refined version of the iteration in \eref{SGS}. At each iteration, the residual term $\LamALn$ is back projected to the scattered data. Consequently, the scattered data are corrected before estimating the scattering source, leading to an improvement in the final solution and stabilization of the iterative process. 
The iteration presented in  \eref{Pair5_iter} can be equivalent to the multiplier-oriented AL iteration, \eref{MOALM}. This equivalence arises for the case $\alpha=1/\regp$ and when we solve the multiplier equation in \eref{Pair5_itera} using a damped LS approach with a damping parameter $\regp$. In this context, the associated Lagrange multipliers are related by the equation $\Lam=(1/\regp)\Lamm$ and $\LamALn=(1/\regp)\LamAL$. 

To further enhance memory efficiency, an alternative approach is to introduce Lagrange multipliers in the data space rather than in the source space. By working in the data space, which typically has a significantly smaller dimension compared to the source space, we can achieve improved memory efficiency without compromising the accuracy of the iterative solution \citep[see][]{Gholami_2022_EFW, Gholami_2023_MWI}.
\section{CONCLUSION}
The quadratic penalty and augmented Lagrangian (AL) methods offer two distinct pathways for regularizing the standard Lagrangian-based full waveform inversion (FWI) process. These methods can be implemented through two equivalent approaches: the first involves an iterative algorithm that computes least-squares (LS) wavefields, projected subsequently to the source space for obtaining multipliers; the second yields another iterative algorithm that calculates LS multipliers, subsequently employed as secondary volume sources to compute the wavefields. Both algorithms update the model by correlating the resultant wavefields and multipliers at each spatial point.

Despite ultimately arriving at the same results, these approaches provide unique perspectives and trade-offs in terms of computational efficiency and complexity. The multiplier-oriented approach offers specific advantages: (1) extending conventional FWI procedures, (2) applicable in both frequency and time domains, and (3) affording a tangible physical interpretation.

We have shown that Lagrange multipliers can be viewed as auxiliary variables, presenting scattering terms inserted into FWI to split and linearize the governing nonlinear equations. Therefore, all of the efforts are focused on creating reliable techniques for multiplier estimate that are accurate during the iterative process.
This interpretation highlights the essential role that multipliers play in improving each method's algorithmic performance and makes it easier to compare them. This study close the gap between theoretical understanding and practical applicability by clarifying the physical relevance of Lagrange multipliers and their implications. Also, it presents useful ways to improve the functionality and applicability of FWI algorithms.

\appendix
\section{FWI method}
In the standard reduced space FWI, the wavefield in obtained by solving the wave equation, $\u(\bold{m})= \bold{A(m)}^{-1}\bold{b}$ and the multiplier is then obtained by solving the adjoint equation, $\Lamm(\bold{m}) = \bold{A(m)}^{-T}\bold{P}^T[\bold{d}-\bold{P}\u(\bold{m})]$.
Due to this strategy, the Lagrangian min-max problem \eref{L} is reduced to a pure minimization problem involving only $\bold{m}$; hence, the term reduced form.
\begin{equation} \label{FWI_obj}
\min_{\bold{m}} \frac{1}{2}\|\bold{P}\bold{A(m)}^{-1}\bold{b}-\bold{d}\|_2^2.
\end{equation}
The gradient of this reduced objective function with respect to the model parameters is
\begin{equation} \label{FWI_grad}
\left(\frac{\partial \bold{A(m)}}{\partial \bold{m}}\u \right)^T \bold{A(m)}^{-T}\bold{P}^T(\bold{d}-\bold{P}\u),
\end{equation}
where $\u=\bold{A(m)}^{-1}\bold{b}$.

\section{WRI method}
The quadratic penalty formulation of the constrained problem in \eref{EFWI_obj}  is
\begin{equation} \label{Penalty_obj}
\min_{\bold{m,u}} \frac{1}{2}\|\bold{Pu-d}\|_2^2+\frac{\regp}{2}\|\bold{A}(\bold{m})\up-\bold{b}\|_2^2.
\end{equation}
where $\regp$ is the penalty parameter. The WRI solves this objective function by alternating between minimization for $\bold{u}$ and $\bold{m}$, which leads to the following iteration \citep{VanLeeuwen_2013_MLM}:
\begin{subequations} \label{WRI}
\begin{align} 
&\u^+ = \arg\min_{\bold{u}} \frac{1}{2}\|\bold{Pu-d^*}\|_2^2+\frac{\regp}{2}\|\bold{A}(\bold{m})\up-\bold{b}^*\|_2^2,\label{WRI_a}\\
&\bold{m}^+ = \arg\min_{\bold{m}} \frac{\regp}{2}\|\bold{A}(\bold{m})\up^+-\bold{b}^*\|_2^2. \label{WRI_b}
\end{align}
\end{subequations} 

\section{IR-WRI method}
The augmented Lagrangian formulation of the problem \eref{EFWI_obj} is 
\begin{equation} \label{AL_obj}
\min_{\bold{m,u}}\max_{\LamAL} \frac{1}{2}\|\bold{Pu-d}\|_2^2+\frac{\regp}{2}\|\bold{A}(\bold{m})\up-\bold{b}\|_2^2 +\langle \LamAL,\bold{A(m)}\u-\bold{b}\rangle.
\end{equation}
Solving this objective function by the alternating direction method of multipliers gives the IR-WRI \citep{Aghamiry_2019_IWR}
\begin{subequations} \label{IRWRI}
\begin{align} 
&\u^+ = \arg\min_{\bold{u}} \frac{1}{2}\|\bold{Pu-d^*}\|_2^2+\frac{\regp}{2}\|\bold{A}(\bold{m})\up-\bold{b}^*\|_2^2+\langle \LamAL,\bold{A}(\bold{m})\up-\bold{b}^*\rangle,\label{IRWRI_w}\\
&\bold{m}^+ = \arg\min_{\bold{m}} \frac{\regp}{2}\|\bold{A}(\bold{m})\up^+-\bold{b}\|_2^2+\langle \LamAL,\bold{A}(\bold{m})\up-\bold{b}^*\rangle, \label{IRWRI_m}\\
&\LamAL^+ = \LamAL+ \regp(\bold{A}(\bold{m}^+)\up^+-\bold{b}). \label{IRWRI_lambda}
\end{align}
\end{subequations}


\begin{thebibliography}{41}
\expandafter\ifx\csname natexlab\endcsname\relax\def\natexlab#1{#1}\fi

\bibitem[Aghamiry et~al.(2019)Aghamiry, Gholami, \& Operto]{Aghamiry_2019_IWR}
Aghamiry, H., Gholami, A., \& Operto, S., 2019.
\newblock Improving full-waveform inversion by wavefield reconstruction with
  alternating direction method of multipliers, {\it Geophysics\/}, {\bf 84(1)},
  R139--R162.

\bibitem[{\SortNoop{Akcelik}}Ak\c{c}elik(2002)]{Akcelik_2002_MNK}
{\SortNoop{Akcelik}}Ak\c{c}elik, V., 2002.
\newblock {\it Multiscale {N}ewton-{K}rylov methods for inverse acoustic wave
  propagation\/}, Ph.D. thesis, Carnegie Mellon University, Pittsburgh,
  Pennsylvania.

\bibitem[Alkhalifah \& Song(2019)]{Alkhalifah_2019_AEW}
Alkhalifah, T. \& Song, C., 2019.
\newblock An efficient wavefield inversion: Using a modified source function in
  the wave equation, {\it Geophysics\/}, {\bf 84}(6), R909--R922.

\bibitem[Barnier et~al.(2023)Barnier, Biondi, Clapp, \&
  Biondi]{Barnier_2023_FWIt}
Barnier, G., Biondi, E., Clapp, R.~G., \& Biondi, B., 2023.
\newblock Full waveform inversion by model extension: theory, design and
  optimization, {\it Geophysics\/}, {\bf 88}(5), 1--206.

\bibitem[Gauthier et~al.(1986{\natexlab{a}})Gauthier, Virieux, \&
  Tarantola]{Gauthier_1986_TDN}
Gauthier, O., Virieux, J., \& Tarantola, A., 1986{\natexlab{a}}.
\newblock Two-dimensional nonlinear inversion of seismic waveforms: numerical
  results, {\it Geophysics\/}, {\bf 51}(7), 1387--1403.

\bibitem[Gauthier et~al.(1986{\natexlab{b}})Gauthier, Virieux, \&
  Tarantola]{Gauthier_1986_TNI}
Gauthier, O., Virieux, J., \& Tarantola, A., 1986{\natexlab{b}}.
\newblock Two-dimensional nonlinear inversion of seismic waveform : numerical
  results, {\it Geophysics\/}, {\bf 51}, 1387--1403.

\bibitem[Gholami et~al.(2022{\natexlab{a}})Gholami, Aghamiry, \&
  Operto]{Gholami_2022_EFW}
Gholami, A., Aghamiry, H.~S., \& Operto, S., 2022{\natexlab{a}}.
\newblock Extended full waveform inversion in the time domain by the augmented
  {L}agrangian method, {\it Geophysics\/}, {\bf 87}(1), R63--R77.

\bibitem[Gholami et~al.(2022{\natexlab{b}})Gholami, Aghamiry, \&
  Operto]{Gholami_2022_OCN}
Gholami, A., Aghamiry, H.~S., \& Operto, S., 2022{\natexlab{b}}.
\newblock On the connection between {WRI} and {FWI}: Analysis of the nonlinear
  term in the hessian matrix, in {\em SEG Technical Program Expanded Abstracts
  2022\/}, pp. 1--5, Society of Exploration Geophysicists.

\bibitem[Gholami et~al.(2023)Gholami, Aghamiry, \& Operto]{Gholami_2023_MWI}
Gholami, A., Aghamiry, H.~S., \& Operto, S., 2023.
\newblock Multiplier waveform inversion ({MWI}): A reduced-space {FWI} by the
  method of multipliers, {\it Geophysics\/}, {\bf 88}(3), R339–R354.

\bibitem[Gill \& Robinson(2012)]{Gill_2012_PDA}
Gill, P.~E. \& Robinson, D.~P., 2012.
\newblock A primal-dual augmented lagrangian, {\it Computational Optimization
  and Applications\/}, {\bf 51}(1), 1--25.

\bibitem[Golub \& Pereyra(2003)]{Golub_2013_VPM}
Golub, G. \& Pereyra, V., 2003.
\newblock Separable nonlinear least squares: the variable projection method and
  its applications, {\it Inverse problems\/}, {\bf 19}(2), R1.

\bibitem[Guttman(1946)]{Guttman_1946_EMC}
Guttman, L., 1946.
\newblock Enlargement methods for computing the inverse matrix, {\it The annals
  of mathematical statistics\/}, pp. 336--343.

\bibitem[Haber et~al.(2000)Haber, Ascher, \& Oldenburg]{Haber_2000_OTS}
Haber, E., Ascher, U.~M., \& Oldenburg, D., 2000.
\newblock On optimization techniques for solving nonlinear inverse problems,
  {\it Inverse problems\/}, {\bf 16}(5), 1263.

\bibitem[Hoffmann et~al.(2021)Hoffmann, Monteiller, \&
  Bellis]{Hoffmann_2021_PFA}
Hoffmann, A., Monteiller, V., \& Bellis, C., 2021.
\newblock A penalty-free approach to {PDE} constrained optimization:
  application to an inverse wave problem, {\it Inverse Problems\/}, {\bf
  37}(5), 055002.

\bibitem[Huang et~al.(2018)Huang, Nammour, \& Symes]{Huang_2018_VSE}
Huang, G., Nammour, R., \& Symes, W.~W., 2018.
\newblock Volume source-based extended waveform inversion, {\it Geophysics\/},
  {\bf 83}(5), R369--387.

\bibitem[Li \& Demanet(2016)]{Li_2016_FWI}
Li, Y.~E. \& Demanet, L., 2016.
\newblock Full-waveform inversion with extrapolated low-frequency data, {\it
  Geophysics\/}, {\bf 81}(6), R339--R348.

\bibitem[Lin et~al.(2023)Lin, van Leeuwen, Liu, Sun, \& Xing]{Lin_2023_FWR}
Lin, Y., van Leeuwen, T., Liu, H., Sun, J., \& Xing, L., 2023.
\newblock A fast wavefield reconstruction inversion solution in the frequency
  domain, {\it Geophysics\/}, {\bf 88}(3), R257–R267.

\bibitem[Lippmann \& Schwinger(1950)]{Lippmann_1950_VPS}
Lippmann, B.~A. \& Schwinger, J., 1950.
\newblock Variational principles for scattering processes. {I}, {\it Physical
  Review\/}, {\bf 79}(3), 469.

\bibitem[M{\'e}tivier et~al.(2017)M{\'e}tivier, Brossier, Operto, \&
  J.]{Metivier_2017_TRU}
M{\'e}tivier, L., Brossier, R., Operto, S., \& J., V., 2017.
\newblock Full waveform inversion and the truncated {N}ewton method, {\it SIAM
  Review\/}, {\bf 59}(1), 153--195.

\bibitem[M{\'e}tivier et~al.(2022)M{\'e}tivier, Brossier, Kpadonou, Messud, \&
  Pladys]{Metivier_2022_ARO}
M{\'e}tivier, L., Brossier, R., Kpadonou, F., Messud, J., \& Pladys, A., 2022.
\newblock A review of the use of optimal transport distances for high
  resolution seismic imaging based on the full waveform, {\it arXiv preprint
  arXiv:2204.08514\/}.

\bibitem[Nocedal \& Wright(2006)]{Nocedal_2006_NO}
Nocedal, J. \& Wright, S.~J., 2006.
\newblock {\it Numerical Optimization\/}, Springer, 2nd edn.

\bibitem[Operto et~al.(2023)Operto, Gholami, Aghamiry, Guo, Mamfoumbi, Beller,
  Aghazade, Mamfoumbi, Combe, \& Ribodetti]{Operto_2023_FWI}
Operto, S., Gholami, A., Aghamiry, H.~S., Guo, G., Mamfoumbi, F., Beller, S.,
  Aghazade, K., Mamfoumbi, F., Combe, L., \& Ribodetti, A., 2023.
\newblock Extending the search space of full-waveform inversion beyond the
  single-scattering born approximation: A tutorial review, {\it Geophysics\/},
  {\bf 88}(6), 1--32.

\bibitem[Plessix(2006)]{Plessix_2006_RAS}
Plessix, R.~E., 2006.
\newblock A review of the adjoint-state method for computing the gradient of a
  functional with geophysical applications, {\it Geophysical Journal
  International\/}, {\bf 167}(2), 495--503.

\bibitem[Pratt et~al.(1998)Pratt, Shin, \& Hicks]{Pratt_1998_GNF}
Pratt, R.~G., Shin, C., \& Hicks, G.~J., 1998.
\newblock {G}auss-{N}ewton and full {N}ewton methods in frequency-space seismic
  waveform inversion, {\it Geophysical Journal International\/}, {\bf 133},
  341--362.

\bibitem[Rizzuti et~al.(2021)Rizzuti, Louboutin, Wang, \&
  Herrmann]{Rizzuti_2021_DFW}
Rizzuti, G., Louboutin, M., Wang, R., \& Herrmann, F.~J., 2021.
\newblock A dual formulation of wavefield reconstruction inversion for
  large-scale seismic inversion, {\it Geophysics\/}, {\bf 86}(6), R879--R893.

\bibitem[Rockafellar(1976)]{Rockafellar_1976_MOA}
Rockafellar, R.~T., 1976.
\newblock Monotone operators and the proximal point algorithm, {\it SIAM
  journal on control and optimization\/}, {\bf 14}(5), 877--898.

\bibitem[Schuster(1985)]{Schuster_1985_HBB}
Schuster, G.~T., 1985.
\newblock A hybrid {BIE+ Born} series modeling scheme: Generalized born series,
  {\it The Journal of the Acoustical Society of America\/}, {\bf 77}(3),
  865--879.

\bibitem[Symes(2020)]{Symes_2020_WRI}
Symes, W.~W., 2020.
\newblock Wavefield reconstruction inversion: an example, {\it Inverse
  Problems\/}, {\bf 36}(10), 105010.

\bibitem[Tarantola(1984)]{Tarantola_1984_ISR}
Tarantola, A., 1984.
\newblock Inversion of seismic reflection data in the acoustic approximation,
  {\it Geophysics\/}, {\bf 49}(8), 1259--1266.

\bibitem[Tarantola(1986)]{Tarantola_1986_SNL}
Tarantola, A., 1986.
\newblock A strategy for non linear inversion of seismic reflection data, {\it
  Geophysics\/}, {\bf 51}(10), 1893--1903.

\bibitem[Tarantola(1988)]{Tarantola_1988_TBI}
Tarantola, A., 1988.
\newblock Theoretical background for the inversion of seismic waveforms
  including elasticity and attenuation, {\it Pure and Applied Geophysics\/},
  {\bf 128}, 365--399.

\bibitem[Tarantola(2005)]{Tarantola_2005_IPT}
Tarantola, A., 2005.
\newblock {\it Inverse {P}roblem {T}heory and {M}ethods for {M}odel {P}arameter
  {E}stimation\/}, Society for {I}ndustrial and {A}pplied {M}athematics,
  Philadelphia.

\bibitem[{van den Berg} \& Kleinman(1997)]{vandenBerg_1997_CSI}
{van den Berg}, P.~M. \& Kleinman, R.~E., 1997.
\newblock A contrast source inversion method, {\it Inverse Problems\/}, {\bf
  13}(6), 1607.

\bibitem[van Leeuwen(2019)]{vanLeeuwen_2019_ANO}
van Leeuwen, T., 2019.
\newblock A note on extended full waveform inversion, {\it arXiv preprint
  arXiv:1904.00363\/}.

\bibitem[{van Leeuwen} \& Herrmann(2016)]{vanLeeuwen_2016_PMP}
{van Leeuwen}, T. \& Herrmann, F., 2016.
\newblock A penalty method for {PDE}-constrained optimization in inverse
  problems, {\it Inverse Problems\/}, {\bf 32(1)}, 1--26.

\bibitem[{van Leeuwen} \& Herrmann(2013)]{VanLeeuwen_2013_MLM}
{van Leeuwen}, T. \& Herrmann, F.~J., 2013.
\newblock Mitigating local minima in full-waveform inversion by expanding the
  search space, {\it Geophysical Journal International\/}, {\bf 195(1)},
  661--667.

\bibitem[Virieux \& Operto(2009)]{Virieux_2009_OFW}
Virieux, J. \& Operto, S., 2009.
\newblock An overview of full waveform inversion in exploration geophysics,
  {\it Geophysics\/}, {\bf 74}(6), WCC1--WCC26.

\bibitem[Wang et~al.(2016)Wang, Yingst, Farmer, \& Leveille]{Wang_2016_FIR}
Wang, C., Yingst, D., Farmer, P., \& Leveille, J., 2016.
\newblock Full-waveform inversion with the reconstructed wavefield method, in
  {\em SEG Technical Program Expanded Abstracts\/}, pp. 1237--1241.

\bibitem[Warner \& Guasch(2016)]{Warner_2016_AWI}
Warner, M. \& Guasch, L., 2016.
\newblock Adaptive waveform inversion: Theory, {\it Geophysics\/}, {\bf 81}(6),
  R429--R445.

\bibitem[Wright \& Recht(2022)]{Wright_2022_ODA}
Wright, S.~J. \& Recht, B., 2022.
\newblock {\it Optimization for data analysis\/}, Cambridge University Press.

\bibitem[Yang et~al.(2018)Yang, Engquist, Sun, \& Hamfeldt]{Yang_2017_AOT}
Yang, Y., Engquist, B., Sun, J., \& Hamfeldt, B.~F., 2018.
\newblock Application of optimal transport and the quadratic {W}asserstein
  metric to full-waveform inversion, {\it Geophysics\/}, {\bf 83}(1), R43--R62.

\end{thebibliography}
\newcommand{\SortNoop}[1]{}

\end{document}